\documentclass[aip]{revtex4-2}

\usepackage[svgnames,rgb]{xcolor} 

\usepackage{comment} 
\usepackage[justification=centering,textfont=it]{caption}
\usepackage{subcaption}
\usepackage{dcolumn}

\usepackage{amsthm}
\usepackage{amsmath}
\usepackage{amssymb}
\usepackage{physics}
\usepackage{bm}
\usepackage{algorithm,algorithmicx}
\usepackage{pgfplots}
    \pgfplotsset{compat=1.14}
\usepackage{tikz}
    \usetikzlibrary{fit,matrix,positioning, decorations.pathreplacing,calc,shapes,arrows,shadows,patterns}

\usepackage[noend]{algpseudocode}
\usepackage{graphicx}

\usepackage{hyperref}
\hypersetup{
	colorlinks=true, 
	breaklinks=true, 
	urlcolor= black, 
	linkcolor= black, 
	citecolor= black, 
	pdftitle={Two methods to approximate the Koopman operator with a reservoir computer}, 
	pdfauthor={Marvyn Gulina, Alexandre Mauroy}, 
	pdfsubject={Simulation} 
	}

\newtheorem{remark}{Remark}

\newcommand{\R}{\mathbb{R}}
\newcommand{\Win}{\bm{W}_{\text{in}}}

\newcommand{\Wout}{\bm{W}_{\text{out}}}
\newcommand{\Wdic}{\bm{W}_1}
\newcommand{\Wcon}{\bm{W}_2}

\newcommand{\bvec}[1]{\bm{#1}}
\newcommand{\bmat}[1]{\bm{#1}}

\newcommand{\C}{\mathbb{C}}

\renewcommand{\sf}{\mathcal{F}}
\newcommand{\sk}{\mathcal{K}}

\newcommand{\sP}{\mathcal{P}}
\newcommand{\sx}{\mathcal{X}}

\makeatletter
\newenvironment{breakablealgorithm}
{
	\begin{center}
		\refstepcounter{algorithm}
		\hrule height.8pt depth0pt \kern2pt
		\renewcommand{\caption}[2][\relax]{
			{\raggedright\textbf{\ALG@name~\thealgorithm} ##2\par}%
			\ifx\relax##1\relax 
			\addcontentsline{loa}{algorithm}{\protect\numberline{\thealgorithm}##2}%
			\else 
			\addcontentsline{loa}{algorithm}{\protect\numberline{\thealgorithm}##1}%
			\fi
			\kern2pt\hrule\kern2pt
		}
	}{
		\kern2pt\hrule\relax
	\end{center}
}
\makeatother

\begin{document}
	
	
	\title{Two methods to approximate the Koopman operator with a reservoir computer} 
	
	
	\author{Marvyn Gulina}
	\email[Electronic mail: ]{marvyn.gulina@unamur.be}

	\author{Alexandre Mauroy}
	\email[Electronic mail: ]{alexandre.mauroy@unamur.be}
	
	\affiliation{Department of Mathematics and Namur Institute for Complex Systems (naXys), University of Namur, Belgium}
	
	
	\date{\today}

	\begin{abstract}
	
	    The Koopman operator provides a powerful framework for data-driven analysis of dynamical systems. 
        In the last few years, a wealth of numerical methods providing finite-dimensional approximations of the operator have been proposed (e.g. extended dynamic mode decomposition (EDMD) and its variants).
        While convergence results for EDMD require an infinite number of dictionary elements, recent studies have shown that only few dictionary elements can yield an efficient approximation of the Koopman operator, provided that they are well-chosen through a proper training process. However, this training process typically relies on nonlinear optimization techniques.
        
        In this paper, we propose two novel methods based on a reservoir computer to train the dictionary. These methods rely solely on linear convex optimization.
        We illustrate the efficiency of the method with several numerical examples in the context of data reconstruction, prediction, and computation of the Koopman operator spectrum.
        These results pave the way to the use of the reservoir computer in the Koopman operator framework.\\
        
        \noindent{Keywords: Koopman operator, dictionary learning, reservoir computing, nonlinear dynamics, dynamic mode decomposition.}
        
	\end{abstract}
	
	\maketitle 
	
	\begin{quotation}
		
		The so-called Koopman operator offers the possibility to turn nonlinear dynamical systems into linear ones.
		In this framework, dynamical systems can be studied with systematic linear techniques and, in particular, they are amenable to spectral analysis.
		However there is a price to pay.
		The Koopman operator is infinite-dimensional and must be approximated by a finite-rank operator (i.e. a matrix) as soon as numerical methods come into play.
		This approximation requires to chose a finite-dimensional subspace, a choice which is not necessarily appropriate since it is made \textit{a priori}.
		Recent methods have been proposed, using neural networks to ``learn'' the best finite-dimensional approximation subspace.
		The main drawback of these methods is that they rely on nonlinear optimization.
		In this paper, we propose to obtain a finite-dimensional approximation of the Koopman operator by using a reservoir computer.
		The reservoir computer is a specific recurrent neural network where only the weights of the nodes on the output layer are trained with the data, a training which can be performed with linear, convex optimization.
		Considering either the internal nodes or the output nodes of the reservoir computer to obtain the finite-dimensional approximation subspace, we derive two novel methods that compute a finite-dimensional approximation of the Koopman operator.
		
	\end{quotation}
	
	
	
	%
	%
	
	%
	
	\section{Introduction}
	
	    Dynamical systems theory plays an important role in the context of data analysis.
	    In fact time-series can often be assumed to be generated by an underlying dynamical system, and to be related to the system orbits through a given observation map.
	    In contrast to this classical description, there exists an alternative description in terms of the observation maps themselves, also called \textit{observables}.
	    The dynamics, and in particular the time evolution of the observables, are then described through the so-called Koopman operator, which is a linear (but infinite-dimensional) operator.
	    In this linear setting, it is natural to study the spectral properties of the operator and relate them to the dynamics of the nonlinear system \cite{Mezic2005}.
	    The related notion of Koopman modes decomposition is also useful to study the systems in many contexts (e.g. fluids dynamics \cite{Rowley2009}, power grids \cite{Susuki2011}, epidemiology \cite{Proctor2015}, control \cite{mauroy2020koopman}).
	    
	    A noticeable fact is that the Koopman operator description is conducive to data analysis.
	    In particular, there exist numerical techniques that can be used to compute a finite-dimensional approximation of the Koopman operator from data.
	    Combined with the spectral analysis relying on Koopman modes, these data-driven techniques lead to the so-called (Extended) Dynamic Mode Decomposition ((E)DMD) method \cite{Schmid2010,H_Tu_2014, Williams2015, Drmac2017}.
	    In practice, EDMD techniques require to choose a specific approximation subspace, or equivalently a finite set of \textit{dictionary functions}.
	    This choice is crucial, but has to be made \textit{a priori} and is therefore not necessarily relevant to provide the best approximation of the operator.
	    
	    Recently, \textit{Dictionary learning} methods based on neural networks have been proposed to provide a relevant set of dictionary functions that are trained with the data and yield appropriate finite-dimensional representations of the Koopman operator \cite{Kevrekidis2017,Takeishi2017}.
	    Subsequent developments have also been made in the context of deep learning \cite{Lusch2018,Yeung2019, Pan2020, dogra2020optimizing}.
	    All these learning methods showed good performances, thereby demonstrating the effectiveness of dictionary learning in EDMD methods.
	    However these techniques require nonlinear optimization of the neural networks, while the classical EDMD method merely relies on linear optimization (i.e. linear least squares regression).
	    
	    In this work, we propose a novel dictionary learning method for EDMD, which relies solely on linear optimization.
	    Our key idea is to combine the EDMD method with a reservoir computer \cite{jaeger2001echo}.
	    In more general contexts, reservoir computing compete with other algorithms on hard tasks such as channel equalization\cite{Skibinsky-Gitlin2018}, phoneme recognition\cite{Triefenbach2010}, and prediction \cite{Pathak2018, gulina2018}, amongst others (see the  survey \cite{lukovsevivcius2009survey,lukovsevivcius2012reservoir}).
	    To our knowledge, we propose the first use of a reservoir computer for dictionary learning in the Koopman operator framework.
	    Note that, very recently, the work \cite{bollt2020} has emphasized a connection between the Koopman operator and the reservoir computer, though in a slightly different setting where the reservoir activation function is linear. This reinforces our claim that reservoir computing is relevant in the context of the Koopman operator.
	    Interestingly the recurrent neural network characterizing the reservoir allows to train the dictionary with a dynamical network rather than with a static one.
	    Hence, generated dictionary functions are nonlinear functions of time-delay coordinates, which are particularly relevant for time-delay systems and also bear some similarity to previous Koopman mode decomposition methods based on delayed coordinates (e.g. prony method \cite{SusukiMezic2015}, Hankel DMD \cite{ArabiMezic2017}).
	    We derive two numerical schemes, where the dictionary functions are respectively the internal states and the outputs of the reservoir computer.
	    We illustrate these two methods with several examples in the context of data reconstruction, prediction, and computation of the Koopman operator spectrum.

    	The rest of the paper is organized as follows.
    	Section \ref{sec:preliminaries} provides an introduction to the Koopman opertor framework, the EDMD method, and the reservoir computer.
    	Section \ref{sec:methods} presents our two methods obtained by combining the EDMD method with the reservoir computer.
    	These two methods are illustrated in Section \ref{sec:results} with numerical examples.
    	Finally, concluding remarks are given in Section \ref{sec:Conclusions}.
    	
	\section{\label{sec:preliminaries}Preliminaries}
    	\subsection{The Koopman operator framework} \label{subsec:Koopman}
        	The Koopman operator provides an alternative framework to describe the evolution of nonlinear systems in a purely linear fashion.
        	Consider an autonomous dynamical system 
        	\begin{equation}
            	\label{eq:dyn_syst}
            	\bvec{x}(t+1) = \bvec{F}(\bvec{x}(t)) \qquad \bvec{x} \in \sx \,,
        	\end{equation}
        	where $\sx$ is (an invariant subset of) the state space and $\bvec{F} : \sx \to \sx$ is a nonlinear map.
        	The Koopman operator is defined as the composition \cite{Koopman1931}
        	\begin{equation}\label{eq:KoopmanDef}
        	    \sk f = f \circ \bvec{F}
        	\end{equation}
        	where $f : \sx \to \C$ is an \textit{observable} that belongs some function space $\sf$.
        	In the following, we will assume that $\sf=L^2(\sx)$ and that $\sx$ is a compact set.
        	It is clear from \eqref{eq:KoopmanDef} that the Koopman operator is linear.
        	Also, while \eqref{eq:dyn_syst} describes the nonlinear dynamics of the state in the space $\sx$, \eqref{eq:KoopmanDef} equivalently describes the linear dynamics of the observables in $\sf$.
        	Roughly speaking the system described in the space $\sx$ is \textit{lifted} into the space $\sf$ when it is described in the Koopman's framework.
        	
        	The Koopman operator description can be used for several purposes. For instance, it can be used for prediction. Indeed, the (vector-valued) identity function $\bvec{g}(\bvec{x}) = \bvec{x}$, also called projection maps, is characterized by the linear evolution
        	\begin{equation}\label{eq:KoopmanEvolution}
        	    \bvec{x}(t+1) = \sk \bvec{g}(\bvec{x}(t)) \,.
        	\end{equation}
        	Provided that the Koopman operator associated with the system is known, \eqref{eq:KoopmanEvolution} allows to predict future trajectories. Moreover the spectral properties of the Koopman operator --- namely the eigenvalues $\lambda \in \mathbb{C}$ and the associated eigenfunctions $\phi \in \mathcal{F}$ satisfying $\sk \phi =  \lambda \phi$ --- provide meaningful information on the underlying dynamical system \cite{Mezic2005, Mauroy2013}.
        	In particular, the eigenvalues are related to internal frequencies of the dynamics and the eigenfunctions reveal geometric properties in the state space.
        	These spectral properties can also be used for control \cite{kaiser2017, Huang2018, Korda2018}, stability analysis \cite{Susuki2014, Mauroy2016GlobalStability}, time-series classification \cite{Surana2018}, analysis and training of neural networks \cite{dogra2020optimizing, Manojlovic2020}, and network identification \cite{Mauroy2017}, to list a few.

        	

    	\subsection{\label{subsec:EDMD}	Finite-dimensional approximation of the Koopman operator}
        	Since the Koopman operator is infinite-dimensional, it is natural and often necessary to compute a finite-dimensional approximation.
        	This approximation is given by the so-called \textit{Koopman matrix} $\bmat{K}$, which represents the projection $\sP$ of the operator onto a subspace $\sf_D$ spanned by the basis functions $\psi_k \in \sf$, $k = 1, \dots, D$, also called \textit{dictionary}.
        	More precisely, the $i^{\text{th}}$ row of $\bmat{K}$ is the coordinate vector of $\sP \sk \psi_i$ in the dictionary.
        	If one denotes $\bmat{\psi}(\bvec{x}) = \left(\psi_1(\bvec{x}), \cdots, \psi_D(\bvec{x}) \right)^T$ and $\sk \bmat{\psi}(\bvec{x}) = \left(\sk \psi_1(\bvec{x}), \cdots, \sk \psi_D(\bvec{x}) \right)^T$, one has
        	\begin{equation*}
        	    \sk \bmat{\psi}(\bvec{x}) \approx \sP \sk \bmat{\psi}(\bvec{x}) = \bmat{K} \bmat{\psi}(\bvec{x})\,,
        	\end{equation*}
        	so that one can obtain an approximation of the evolution of the dictionary functions under the action of the Koopman operator.
        	In particular, if the identity belongs to the dictionary, it follows from \eqref{eq:KoopmanEvolution} that an approximation of the system trajectories can be computed.
        	
        	The finite-dimensional approximation of the Koopman operator can be obtained from data through the so-called Extended Dynamic Mode Decomposition (EDMD) method \cite{Williams2015}.
        	Given a set of snapshot pairs $\left\{ (\bvec{x}_t , \bvec{x}'_t = \bvec{F}(\bvec{x}_t)) \right\}_{t = 1}^{T}$, the Koopman matrix is given by
        	\begin{equation}\label{eq:EDMDproblem}
            	\bmat{K} 	= \arg \min_{\tilde{\bmat{K}} \in \R^{D \times D}} \sum_{t=1}^{T} \norm{\bmat{\psi}(\bvec{x}'_t) - \tilde{\bmat{K}}\bmat{\psi}(\bvec{x}_t)}^2
            	= \bmat{\Psi}' \, \bmat{\Psi}^+
        	\end{equation}
        	where ${}^+$ denotes the pseudo-inverse, and where $\bmat{\Psi}$ and $\bmat{\Psi}' \, \in \, \R^{D \times T}$ denote the matrices whose columns are $\bmat{\psi}(\bmat{x}_t)$ and $\bmat{\psi}(\bmat{x}'_t)$, respectively, for $t \in \{1, \cdots, T\}$. The Koopman matrix is the solution to a least squares problem and therefore represents the approximation of the operator obtained with a discrete orthogonal projection.
        	Note that the specific dictionary $\bmat{\psi}(\bvec{x}) = \bvec{x}$  leads to the classical Dynamic Mode Decomposition (DMD) algorithm \cite{Rowley2009,Schmid2010,H_Tu_2014}.
        	This statement justifies the term ``extended'' introduced above.
        	
        	The finite-dimensional approximation of the operator depends on both the projection and the dictionary of basis functions.
        	In a data-driven context, the discrete orthogonal projection used in the EDMD method is a natural and appropriate projection to use.
        	However, the choice of the dictionary is somehow arbitrary but crucial since it affects the quality of the approximation.
        	The original EDMD method\cite{Williams2015} relies on a dictionary that is fixed and chosen \textsl{a priori} (e.g. polynomial functions, radial basis functions).
        	Recently, machine learning techniques have been used to guide the choice of the dictionary\cite{Kevrekidis2017}.
        	Building on this result, we propose to select the dictionary functions through a reservoir computer.
    	
    	\subsection{\label{subsec:RC}Reservoir computer}       	
        	\begin{figure}
        		\centering
        		\includegraphics[width=0.65\columnwidth]{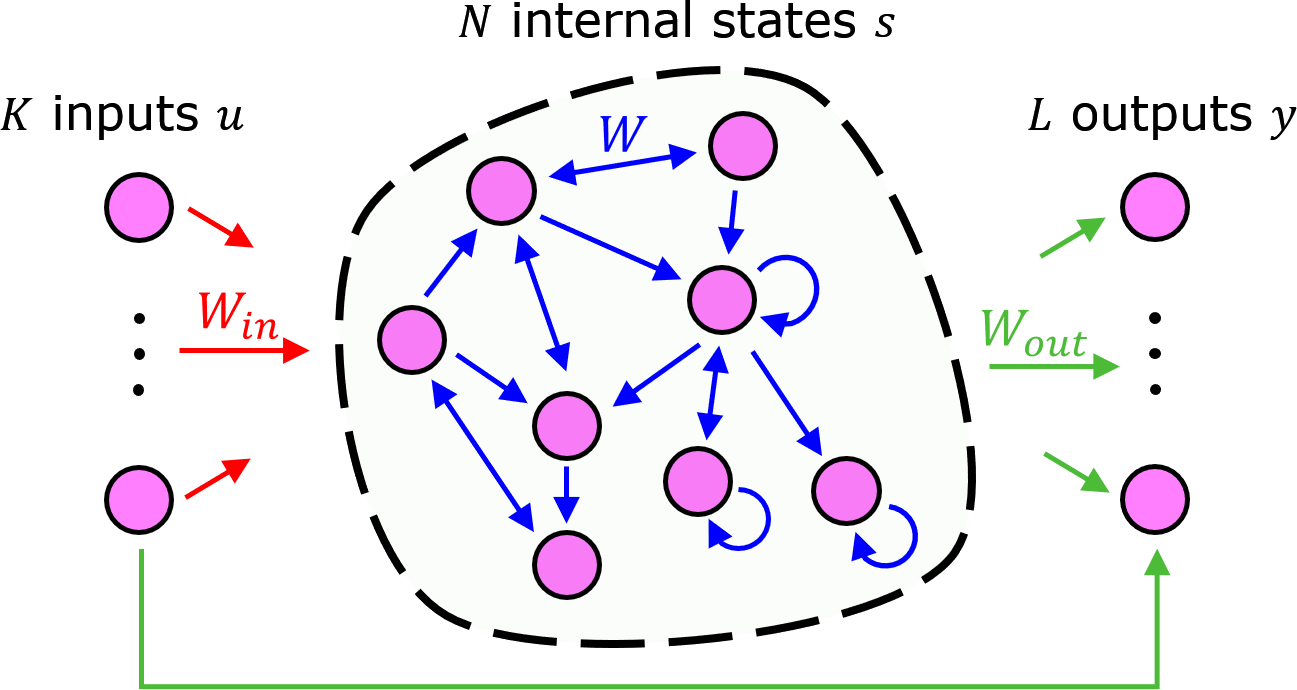}
        		\caption{\label{fig:RC}Layout of the reservoir computer.}
        	\end{figure}
        	
        	The reservoir computer is a discrete-time neural network which consists of three layers: the inputs, the reservoir, and the outputs (FIG. \ref{fig:RC}).
        	We denote the input signals by $\bvec{u} \in \R^K$, the reservoir states by $\bvec{s} \in \R^N$, and the output signals by $\bvec{y}\in \R^L$. 
        	The reservoir states are updated according to the dynamics\cite{jaeger2001echo}
        	\begin{equation}{\label{eq:evolutionDiscretRC}}
            	\bvec{s}(t+1) = (1 - Ca) \bvec{s}(t) + C \tanh \left[\Win  \bvec{u}(t+1) + \bm{W} \bvec{s}(t) + \bvec{\nu}(t) \right]
        	\end{equation}
        	where $C$ is a timescale constant, $a$ is the leaking rate, $\bm{W} \in \R^{N \times N}$ and $\Win \in \R^{N \times K}$ are the matrices of internal connection weights and input weights, respectively, and
            $\bvec{\nu}$ is a noise term.
            The matrix $\bm{W}$ is typically sparse and its density (\textit{i.e.} the proportion of non-zero elements) is denoted by $\gamma$.
        	The nonzero entries of $\Win$ and $\bm{W}$ 
        	are uniformly randomly distributed over $[-1,1]$ and $[-w,w]$, respectively.
        	The components of $\bvec{\nu}$ are uniformly distributed over $[-\varepsilon, \varepsilon]$.
        	This noise term is added according to the work \cite{jaeger2001echo} as an alternative to Tikhonov regularization for the output weights training.
        	
        	The reservoir can contain loops and is therefore a \textit{recurrent neural network}.
        	Furthermore, the gain $w$ is chosen such that the spectral radius $\rho$ of $\bm{W}$ satisfies the echo state property \cite{jaeger2004harnessing,jaeger2001echo}
        	\begin{equation}{\label{eq:ESP}}
        	    |1 - C(a - \rho)| < 1,
        	\end{equation}
        	so that the reservoir ``forgets'' the initial condition $\bvec{s}(0)$, which is uniformly distributed over $[0, 1]$.
        	Hence the computer reservoir is an \textit{echo state network}.
        	In practice, we can discard the first (transient) states corresponding to the initialization of the reservoir.
        	
        	Finally, the outputs are given by
        	\begin{equation}{\label{eq:rcout}}
        	    \bvec{y}(t) = \Wout \bar{\bvec{s}}(t),			
        	\end{equation}
        	where 
        	$\bar{\bvec{s}}(t) =  [\bvec{s}(t) ; \bvec{u}(t)] \in \R^{\bar{N}}$ (with $\bar{N} = N + K$) is the vertical concatenation of internal states and inputs and where $\Wout \in \R^{L \times \bar{N}}$ is the matrix of output weights.
        	It is noticeable that the outputs are obtained through \textit{linear} combinations of the states and that \textit{only} the output weights are trained.
        	This is a computational advantage of the reservoir computer that we will leverage.
	
	\section{\label{sec:methods}Numerical methods}
	
    	
    	In this Section, we present our two methods, which combine the EDMD technique with the reservoir computer.
    	The key idea is to use (linear combinations of) the internal states of the reservoir as dictionary functions.
    	The first method uses all the states of the reservoir, while the second method selects a subset of these states.
    	
    	
    	
    	\subsection{\label{subsec:meth1}Method 1: EDMD using a reservoir computer}
        	A straightforward method consists in using all the reservoir internal states as dictionary functions, \textit{i.e.} $\bmat{\Psi}=\bvec{\bar{s}}$.
        	In this case the optimization problem \eqref{eq:EDMDproblem} becomes
        	\begin{equation}\label{eq:EDMDproblem1}
            	\bmat{K} = \arg \min_{\tilde{\bmat{K}} \in \R^{\bar{N} \times \bar{N}}} \sum_{t=1}^{T-1} \norm{\bvec{\bar{s}}(t+1) - \tilde{\bmat{K}}\bvec{\bar{s}}(t)}^2
        	\end{equation}
        	and its solution is given by
        	\begin{equation}\label{eq:EDMDRC}
        	    \bmat{K} = \bmat{S}' \bmat{S}^{+}
        	\end{equation}
        	where $\bmat{S}$ and $\bmat{S}' \, \in \, \R^{\bar{N} \times (T - 1)}$ denote the matrices whose columns are $\bvec{\bar{s}}(t)$ for $t \in \{1, \cdots, T-1\}$ and $t \in \{2, \cdots, T\}$, respectively.
        	The internal states evolve according to the dynamics \eqref{eq:evolutionDiscretRC}, where the input $\bvec{u}(t)$ is a trajectory $\bvec{x}(t)$ (equivalently, $\bvec{\bar{s}}=[\bvec{s};\bvec{x}]$).
        	For this reason, the data points are generated from a single trajectory, and not from a set of scattered data pairs.
        	
        	\begin{remark} The proposed dictionary can be interpreted as nonlinear functions of time-delay coordinates.
        	Indeed, since the internal states are solutions to \eqref{eq:evolutionDiscretRC}, they depend on past values of the input $\bvec{u}(t)$, or equivalently of the state $\bvec{x}(t)$.
        	Moreover, provided that the reservoir satisfies the echo state property and is initialized for a sufficiently long time before generating the data, the internal states do not depend on the (random) initial condition and are time-independent.
        	The use of time-delay coordinates to construct the basis functions is reminiscent of previous works in the context of the Koopman operator \cite{ArabiMezic2017,SusukiMezic2015,Khodkar2019}.
        	\end{remark}
        	
        	The proposed method is summarized in Algorithm \ref{algo:method1}.
        	
        	\begin{breakablealgorithm}
        		\label{algo:method1}
        		\caption{\label{algo:EDMDRC}EDMD using a reservoir computer}
        		\begin{algorithmic}[1]
        			\Statex{\bf Input:} Sampled trajectory $\bvec{x}(t)$, $t=1,\dots,T$; parameters $N$, $L$, $C$, $a$, $\gamma$, $\rho$, $w_{in}$, $\varepsilon$.
        			\Statex{\bf Output:} Koopman matrix $\bmat{K}$.
        			
        			\State Initialize $\Win$, $\bm{W}$ (see Section \ref{subsec:RC}) and set random initial reservoir states.
        			
        			\State Update the reservoir states according to \eqref{eq:evolutionDiscretRC} with the input $\bvec{u}(t)=\bvec{x}(t)$, for $t=1,\dots,T-1$.
        			
        			\State Compute the Koopman matrix $\bmat{K} = \bmat{S}' \bmat{S}^{+}$ \eqref{eq:EDMDRC} (possibly discarding the first values of the states $\bvec{\bar{s}}$).
        		\end{algorithmic}
        	\end{breakablealgorithm}
    	
    	\subsection{\label{subsec:meth2} Method 2: dictionary learning for EDMD using a reservoir computer}
        	The method presented above yields a $\bar{N} \times \bar{N}$ Koopman matrix, which is not so convenient since the number of internal states is typically large.
        	In order to reduce the size of the Koopman matrix, we propose to define the dictionary functions as the output \eqref{eq:rcout} of the reservoir, that is we select $L \ll N$ linear combinations of the states $\bvec{\bar{s}}=[\bvec{s};\bvec{x}]$ obtained through a dictionary learning step.
        	The optimization problem is written as
        	\begin{equation}\label{eq:kevrekidisProblemWithRC}
            	\min_{\Wout,\bmat{K}}
            	\sum_{t = 1}^{T-1} \norm{\Wout \bvec{\bar{s}}(t+1)  - \bmat{K} \Wout \bvec{\bar{s}}(t)}^2. 
        	\end{equation}
        
        	It is clear that $\Wout=\bmat{0}$ is a trivial solution to the optimization problem.
        	We therefore add the projection maps $\bvec{x}$ to the outputs, considering the augmented output weights matrix $\Wout = [\Wdic ; \Wcon] \in \R^{(L+K) \times \bar{N}}$ with the output weights $\Wdic \in \R^{L \times \bar{N}}$ related to the main dictionary functions and the output weights $\Wcon \in \R^{K \times \bar{N}}$ related to the projection maps $\bvec{x}$.
        	The dictionary size is given by $D = L+K$.
        	Since $\bvec{u}(t)=\bvec{x}(t)$, we have $\Wcon[\bvec{s} ; \bvec{u}] = \bvec{u}$, so that the constraint $\Wcon=[\bmat{0}_{K,N} , \bmat{I}_K]$ (where $\bmat{0}_{K,N}$ is the $K\times N$ zero matrix and $\bmat{I}_K$ is the $K\times K$ identity matrix) prevents the trivial zero solution.
        	In this case, we only need to optimize over the weights $\Wdic$.
        	
        	This method presented can be interpreted as an intermediate method between the first method developed in Section \ref{subsec:meth1} and the DMD method.
        	If we set $L = N$, the dictionary contains the maximal number of independent functions constructed with the internal states of the reservoir.
        	These functions can be chosen as the internal states themselves and we recover the first method.
        	In contrast, in the case $L=0$, there is no optimization performed on the outputs weights and the optimization on $\bmat{K}$ associated with linear basis functions is equivalent to DMD. The proposed method can therefore be seen as a trade-off where an optimal subset of basis functions is obtained through the training process.
        	
        	Similarly to \cite{Kevrekidis2017}, \eqref{eq:kevrekidisProblemWithRC} can be solved with two alternating steps:
        	\begin{enumerate}
        		\item (Computation of the Koopman matrix) fix $\Wdic$ and optimize $\bmat{K}$ ; \label{firstStep}
        		\item (Dictionary learning) fix $\bmat{K}$ and optimize $\Wdic$. \label{secondStep}
        	\end{enumerate}
        	A key advantage is that \textit{both} optimization steps rely on \textit{linear} optimization.
        	It is clear that step \ref{firstStep} is a least squares problem, whose solution is given by
        	\begin{equation}\label{eq:step1Method2}
        	    \bmat{K} = \Wout \bmat{S}' (\Wout \bmat{S})^{+} = \Wout \bmat{S}' \bmat{S}^{+} \Wout^{+} \,.
        	\end{equation}
        	In fact, this is the standard EDMD problem $\bm{K} = \bmat{\Psi}' \, \bmat{\Psi}^+$ in the case of our specific choice of dictionary functions $\bmat{\Psi} = \Wout \bmat{S}$.
        	We can also note that the matrix $\bmat{S}' \bmat{S}^{+}$ in \eqref{eq:step1Method2} is the Koopman matrix computed in the first method.
        	The optimization over the output weight matrix $\Wout$ somehow acts as a coupling in the optimization problem, where the columns of the Koopman matrix are optimized simultaneously in order to minimize the overall residue in \eqref{eq:kevrekidisProblemWithRC}.
        	
        	Step \ref{secondStep} can also be cast into a least squares problem.
        	Denoting
        	\[ 
            	\bmat{K} = \begin{pmatrix}
            	\bmat{K}_{11} & \bmat{K}_{12}\\
            	\bmat{K}_{21} & \bmat{K}_{22}
            	\end{pmatrix}\,,
        	\]
        	we can write problem \eqref{eq:kevrekidisProblemWithRC} as the equation
        	\begin{equation}\label{eq:stepMethod2}
            	\begin{split}
                	\bmat{W} \bmat{S}' - \bmat{K}\bmat{W}\bmat{S} & = 
                	\begin{pmatrix}
                    	\Wdic \bmat{S}'\\
                    	\Wcon \bmat{S}'
                	\end{pmatrix}
                	-
                	\begin{pmatrix}
                    	\bmat{K}_{11}\Wdic \bmat{S} + \bmat{K}_{12}\Wcon \bmat{S}\\
                    	\bmat{K}_{21}\Wdic \bmat{S} + \bmat{K}_{22}\Wcon \bmat{S}
                	\end{pmatrix} \\
                	& = 
                	\begin{pmatrix}
                    	\Wdic \bmat{S}' \\
                    	\bmat{0}_{K,T-1} \end{pmatrix} - \begin{pmatrix}
                    	\bmat{K}_{11} \\
                    	\bmat{K}_{21}
                	\end{pmatrix}
                	\Wdic \bmat{S}
                	\quad - 
                	\begin{pmatrix}
                    	\bmat{K}_{12}\Wcon \bmat{S}\\
                    	- \Wcon \bmat{S}' + \bmat{K}_{22}\Wcon \bmat{S}
                	\end{pmatrix}\\
                	& = \bmat{0}_{D,T-1} \,.
            	\end{split}
        	\end{equation}
        	Denoting the independent terms $\bmat{C}_1=\bmat{K}_{12}\Wcon \bmat{S}$ and $\bmat{C}_2=-\Wcon \bmat{S}' + \bmat{K}_{22}\Wcon \bmat{S}$, we compute the optimal output weights $\Wdic$ as the least squares solution given by
        	\begin{equation}\label{eq:altMethod2}
            	\Wdic(:)
            	= \left[
                	\begin{array}{c}
                    	(\bmat{S}')^T \otimes \bmat{I}_L - \bmat{S}^T \otimes \bmat{K}_{11} \\
                    	-\bmat{S}^T \otimes \bmat{K}_{21}
                	\end{array}
            	\right]^{+}
            	\left[
                	\begin{array}{c}
                    	\bmat{C}_1(:) \\ 
                    	\bmat{C}_2(:)
                	\end{array}
            	\right]\,,
        	\end{equation}
        	where (:) denotes the vectorization of matrices and $\otimes$ is the Kronecker product.
        	
        	A variant of the above numerical scheme can also be obtained. Instead of computing the least squares solution \eqref{eq:altMethod2} in terms of output weights $\bmat{W}_1$, we could replace them by the reservoir outputs $\bmat{\Psi}_{1} = \bmat{W}_1\bmat{S}$ in \eqref{eq:stepMethod2}.
        	Using $\bmat{W}_1 = \bmat{\Psi}_{1} \bmat{S}^{+}$, we obtain the Sylvester equations 
        	\begin{eqnarray*} 
        		\bmat{\Psi}_{1} \bmat{S}^{+} \bmat{S}' - \bmat{K}_{11} \bmat{\Psi}_{1} & = & \bmat{C}_{1} \\
        		-\bmat{K}_{21} \bmat{\Psi}_{1} & = & \bmat{C}_{2} \,,
        	\end{eqnarray*}
        	whose least squares solution is
        	\begin{equation}\label{eq:DLRCsystem}
            	\bmat{\Psi}_1(:)
            	= 	\left[
                	\begin{array}{c}
                    	(\bmat{S}^{+} \bmat{S}')^T \otimes \bmat{I}_{L}
                    	-\bmat{I}_{T-1} \otimes \bmat{K}_{11}  \\
                    	-\bmat{I}_{T-1} \otimes \bmat{K}_{21}
                	\end{array}
            	\right]^{+}
            	\left[
                	\begin{array}{c}
                    	\bmat{C}_1(:) \\ 
                    	\bmat{C}_2(:)
                	\end{array}
            	\right] \,.
        	\end{equation}
        

        	In the case $\bar{N} \ge T-1$, it is noticeable that $\bmat{S}^{+} \bmat{S}'$ is a companion matrix, so that only the last column has to be effectively computed.
        	When $\bar{N} > T-1$, optimizing over $\bmat{\Psi}_{1} = \bmat{W}_1\bmat{S}$ in the variant of the method makes the optimization problem more constrained, since the basis functions are not described anymore as linear combinations of the state reservoirs but through their values at a smaller number $T-1< \bar{N}$ of sample points.
        	Conversely, the case $\bar{N} < T-1$ is a relaxation of the optimization problem which is less computationally efficient since $\bmat{S}^{+} \bmat{S}'$ is not sparse in this case.
        	The intermediate setting $\bar{N} = T-1$, where $\bmat{S}$ is a square matrix and $\bmat{\Psi}_{1} = \bmat{W}_1\bmat{S}$, can be interpreted as a change of variable.
        	It appears as an appropriate choice, but the size of the reservoir becomes forced by the size of the dataset.

        	Both equations \eqref{eq:altMethod2} and \eqref{eq:DLRCsystem} provide an effective way to solve the second step of the alternating optimization, yielding two variants of our proposed second method.
        	The first variant provides the output weights describing the dictionary functions, while the second variant provides the updated values of the dictionary functions, which can directly be used to compute the Koopman matrix $\bmat{K}$.
        	Both variants require to invert a matrix of $D(T-1) \times \bar{N}L$ and $D(T-1) \times (T-1)L$, respectively.
        	This matrix is not sparse for the first variant, but is sparse for the second variant (provided that $\bar{N} \ge T-1$).
        	It follows that the second variant is more efficient from a computational point of view.
   
        	The method with its two variants is summarized in Algorithm \ref{algo:DLRC}.

        	\begin{breakablealgorithm}
        		\caption{\label{algo:DLRC}Dictionary learning for EDMD using a reservoir computer}
        		\begin{algorithmic}[1]
        			\Statex{\bf Input:} Sampled trajectory $\bvec{x}(t)$, $t=1,\dots,T$; parameters $N$, $L$, $C$, $a$, $\gamma$, $\rho$, $w_{in}$, $\varepsilon$.
        			\Statex{\bf Output:} Koopman matrix $\bmat{K}$, output weights $\Wdic$
        			\State Initialize $\Win$, $\bm{W}$ (see Section \ref{subsec:RC}) and set random initial reservoir states and output weights $\Wdic$.
 
        			\State Update the reservoir states according to \eqref{eq:evolutionDiscretRC} with the input $\bvec{u}(t)=\bvec{x}(t)$, for $t=1,\dots,T-1$.
        			
        			\State Compute the dictionary values $\bmat{\Psi} = [\bmat{\Psi}_1;\bmat{\Psi}_2]$ with $\bmat{\Psi}_1 = \Wdic \bmat{S}$ and $\bmat{\Psi}_2=[\bmat{0}_{K,N} , \bmat{I}_K]\bmat{S}$ (possibly discarding the first values of the states $\bvec{\bar{s}}$).		
        			
        			\Loop
            			\State Solve step \ref{firstStep}: compute the Koopman matrix $\bmat{K} = \Wout \bmat{S}' \bmat{S}^{+} \Wout^{+}$ (variant 1) or $\bmat{K} = \bmat{\Psi}' \, \bmat{\Psi}^{+}$ (variant 2).
            			
            			\State Solve step \ref{secondStep}: update the values $\Wdic$ using \eqref{eq:altMethod2} (variant 1) or $\bmat{\Psi}_1$ using \eqref{eq:DLRCsystem} (variant 2).
    
            			\State Stop if some criterion is satisfied (e.g. upper bound on the number of iterations, lower bound on the least squares error).
            			
        			\EndLoop
        			
        			\State If variant 2 is used: compute the output weights $\Wdic = \bmat{\Psi}_1 \, \bmat{S}^+$.
        			
        		\end{algorithmic}
        	\end{breakablealgorithm}

    	\subsection{Application to prediction and spectral properties}   
        	\subsubsection{\label{subsubsec:recoAndPred}Reconstruction and prediction}

            	The Koopman matrix $\bmat{K}$ provided by Algorithm \ref{algo:EDMDRC} or \ref{algo:DLRC} is optimized so that
            	\[
            	    \forall t \in \{1, \dots, T-1 \} : \bmat{\Psi}(t+1) \approx \bmat{K} \bmat{\Psi}(t).
            	\]
            	
                It follows that we can iterate the matrix to recompute known dictionary values $\hat{\bvec{\Psi}}(t+1) = \bmat{K}^t \ \bvec{\Psi}(1)$ (reconstruction) or predict new values  $\hat{\bvec{\Psi}}(t+T) = \bmat{K}^t \ \bvec{\Psi}(T)$ from the last data point (prediction). In particular, predicted states are obtained by considering the values of the dictionary functions related to the projection maps.
                
                It should be noted that our use of the reservoir computer differs from the classic use in the context of prediction and is not aimed at this specific prediction objective.
            	Here the outputs are not optimized so that they provide the best predictions of the state.
            	Instead they provide the basis functions that yield the most accurate approximation of the Koopman operator, which can in turn be used for prediction as a by-product.
            	This observation is particularly relevant to Method 2 and will be discussed with more details in Section \ref{sec:results}.
                
                In order to test the quality of the Koopman matrix approximation, we will first compute the optimization residue
                \begin{equation}
                    \label{eq:residue}
                \sum_{t=1}^{T} \norm{\bmat{\psi}(\bvec{x}'_t) - \tilde{\bmat{K}}\bmat{\psi}(\bvec{x}_t)}^2.
                \end{equation}
                We will also consider the reconstruction error
            	\begin{equation}
            	\label{eq:error}
            	    \bvec{E}(t) = \bmat{K}^t \ \bvec{\Psi}(1) - \bvec{\Psi}(t+1)
            	\end{equation}
            	and in particular the Normalized Root Mean Square Error (NRMSE)
            	\[
            	    \bvec{NRMSE} = \frac{\sqrt{\sum_{t=1}^{T-1} \bvec{E}(t)^2 }}%
            	    {\sqrt{\sum_{t=1}^{T-1} \left[ \bvec{\Psi}(t+1) -  \frac{1}{T-1}  \sum_{\tau=1}^{T-1} \bvec{\Psi}(\tau+1) \right]^2 }}
            	\]
            	where the square operations, the square roots and the quotient are considered element-wise.
            	The NRMSE value can be interpreted as follows: $NRMSE = 0$ means that the two series perfectly match and $NRMSE = 1$ is the error obtained when the reconstructed time serie is a constant value equal to the mean value of the other one.
            	Similarly we will denote by $\bvec{nrmse}$ the error restricted to the projection maps.

        	\subsubsection{Spectral properties}
            	The Koopman matrix $\bmat{K}$ can be used to compute spectral properties of the Koopman operator $\sk$\cite{Williams2015}.
            	In particular, the eigenvalues of $\bmat{K}$ provide an approximation of the Koopman operator spectrum.
            	Its left eigenvectors $\bvec{w}$ provide the expansion of Koopman eigenfunctions $\phi$ in the basis given by the dictionary functions, \textit{i.e.} $\phi \approx \bvec{w}^T \bvec{\psi}$.
            	Note also that the right eigenvectors are related to the Koopman modes.

	\section{\label{sec:results}Results}  
	        
	        In this Section, we illustrate the performance of our methods with several datasets in the context of trajectory reconstruction, prediction, and computation of spectral properties.

        \subsection{\label{subsec:resultsPreliminaries} Datasets and parameters}	        
             \subsubsection{Dynamical systems and data generation}    
             
                We consider several systems, including chaotic dynamics.

                \paragraph{Van der Pol system.}
                Using the limit-cycle dynamics        
                \begin{equation}
                    	\begin{array}{rcl}
                        	\Dot{x}_1 & = & x_2\\
                        	\Dot{x}_2 & = & \mu(1 - x_1^2)x_2 - x_1\\
                    	\end{array}
                	\end{equation}
                	with $\mu = 1$, we have generated $T = 501$ data points over the time interval $[0,20]$ (time step $h = 0.04$) for the initial condition $\bvec{x}(0) = (-4,5)$. \\
            
                \paragraph{Duffing system.}
                The dynamics
                	\begin{equation}
                    	\begin{array}{rcl}
                        	\Dot{x}_1 & = & x_2\\
                        	\Dot{x}_2 & = & - \gamma x_2 - (\alpha x_1^2 + \beta ) x_1\\
                    	\end{array}
                	\end{equation}
                	with $\alpha = 1$, $\beta = -1$, and $\gamma = 0.5$
                	admit a stable equilibrium at the origin.
                	We have generated $T = 501$ data points over the time interval $[0,20]$ (time step $h = 0.04$) for the initial condition $\bvec{x}(0) = (-1.21, 0.81)$.\\
                	
               \paragraph{Mackey-Glass system.}
                We consider the following delayed equation\cite{MackeyGlassSystem}:
                	\begin{equation}
                	    \Dot{x} = \displaystyle \frac{\alpha x_\tau}{(1 + x_\tau^{n})} - \beta x\\
                	\end{equation}
                	where $x_\tau = x(t - \tau)$ with $\tau = 17$ and $\alpha = 0.2$, $\beta = 0.1$, and $n = 10$.
                	We have generated $T = 501$ data points over the time interval $[0,500]$ (time step $h = 1$) for the initial condition $x(t < 0) = 0.1$.
                	Note that the system is integrated using the Matlab function \texttt{dde23}. \\

                \paragraph{R\"{o}ssler system.}
                We have used the chaotic dynamics
                	\begin{equation}
                    	\begin{array}{rcl}
                        	\Dot{x}_1 & = & - x_2 - x_3\\
                        	\Dot{x}_2 & = & x_1 + \alpha x_2\\
                        	\Dot{x}_3 & = & \beta + (x_1 - \gamma) x_3\\
                    	\end{array}
                	\end{equation}
                	with $\alpha = 0.1$, \quad $\beta = 0.1$, \quad $\gamma = 14$ to generate $T = 601$ data points over the time interval $[0,300]$ (time step $h = 0.5$) for the initial condition $\bvec{x}(0) = (2, 1, 5)$. \\
            	
                \paragraph{Lorenz-63 system.}
                Using the chaotic dynamics
                	\begin{equation}
                    	\begin{array}{rcl}
                        	\Dot{x}_1 & = & s(x_2 - x_1)\\
                        	\Dot{x}_2 & = & r x_1 - x_2 - x_1 x_3\\
                        	\Dot{x}_3 & = & x_1 x_2 - b x_3\\
                    	\end{array}
                	\end{equation}
                	with $s = 10$, \quad $r = 28$, \quad $b = 8/3$,
                	we have generated $T = 751$ data points over the time interval $[0,15]$ (time step $h = 0.02$) for the initial condition $\bvec{x}(0) = (3, 3, 19)$.\\
                	
                	\begin{remark}\label{rmq:dataRescaled}
                	    The data are scaled through a linear transformation that maps the minimum and maximum values of each state to $-1$ and $1$, respectively.
                	    The values of the initial conditions are given  before rescaling.
                	    The data that support the findings of this study are available from the corresponding author upon reasonable request.
                	\end{remark}
                	
           \subsubsection{Parameter values}
                	          
                The number of basis functions is set to $D = N + K = 1000$ and $D = L + K = 15$ for the first and the second method, respectively.
                Note that the numbers $N$ and $L$ depend on the state dimension $K$.
                We note that $\tilde{N} = 1000 > T - 1$ for all study cases.
                Although this choice yields a more constrained optimization problem  (for the second variant of the second method, see Section \ref{subsec:meth2}), it yields the best results in terms of reconstruction, prediction and spectral properties while remaining computationally efficient.
                
                For the second method, the second variant is used and the number of iterations is limited to 20.
                In the examples, the EDMD method is also used for comparison purpose with a dictionary of $N$ Gaussian radial basis functions $\psi_k(\bvec{x})=e^{-\gamma \|\bvec{x}-\hat{\bvec{x}}_k\|^2}$ (with $D = N + K = 1000$), where $\hat{\bvec{x}}_k$ is the center and with $\gamma=0.05$.
                
                In most cases, the reservoir parameters are kept constant for every systems.
                The spectral radius of the internal weights is set to $\rho=0.79$ for all cases.
                The leaking rate is set to $a = 3$ for all systems except for the Mackey-Glass system where $a=1$.
                The noise level in the reservoir is $\varepsilon = 10^{-4}$.
                The time constant is set to $C=0.45$ for the Van der Pol system and the Duffing system, and is set to $C=0.11$ for the Mackey-Glass system, the R\"{o}ssler system, and the Lorenz-63 system.

        \subsection{\label{subsec:reco} Reconstruction results}
        
\begin{figure}
                \centering
                \begin{tikzpicture}
                    \begin{axis}[
                        x tick label style={
                        /pgf/number format/1000 sep=},
                        ylabel={Residue},
                        width = 0.7\columnwidth,
                        ymode=log,
                        legend style={
                            at={(0.5,-0.1)},
                            anchor=north,
                            legend columns=-1,
                            font=\scriptsize},
                        symbolic x coords={
                            V.d.P.,
                            Duffing,
                            M.-G.,
                            R\"{o}ssler,
                            Lorenz,
                            Dummy},
                        xtick=data,
                        tick label style={font=\scriptsize},
                        ybar interval=0.5,
                        ]
                        
                        \addplot[fill=black, postaction={pattern=dots, pattern color = white}] coordinates {
                            (V.d.P.,3.18e-06)
                            (Duffing,1.09e-08)
                            (M.-G.,1.56e+00)
                            (R\"{o}ssler,6.51e-06)
                            (Lorenz,7.26e-06)
                            (Dummy,1)};
                        
                        \addplot[fill=red, postaction={pattern=north east lines}] coordinates {
                            (V.d.P.,5.44e-20)
                            (Duffing,8.28e-26)
                            (M.-G.,3.28e-24)
                            (R\"{o}ssler,4.32e-22)
                            (Lorenz,8.56e-21)
                            (Dummy,1)};
                        
                        \addplot[fill=green, postaction={pattern=grid}] coordinates {
                            (V.d.P.,2.08e-14)
                            (Duffing,3.15e-12)
                            (M.-G.,1.42e-12)
                            (R\"{o}ssler,2.09e-13)
                            (Lorenz,7.02e-11)
                            (Dummy,1)};
                        
                        \legend{EDMD (RBF),Method 1,Method 2}
                    \end{axis}
                \end{tikzpicture}
                \caption{\label{fig:residues}
                The optimization residues of \eqref{eq:EDMDproblem} are computed for the different systems and show that our methods 1 and 2 outperform EDMD.}
            \end{figure}
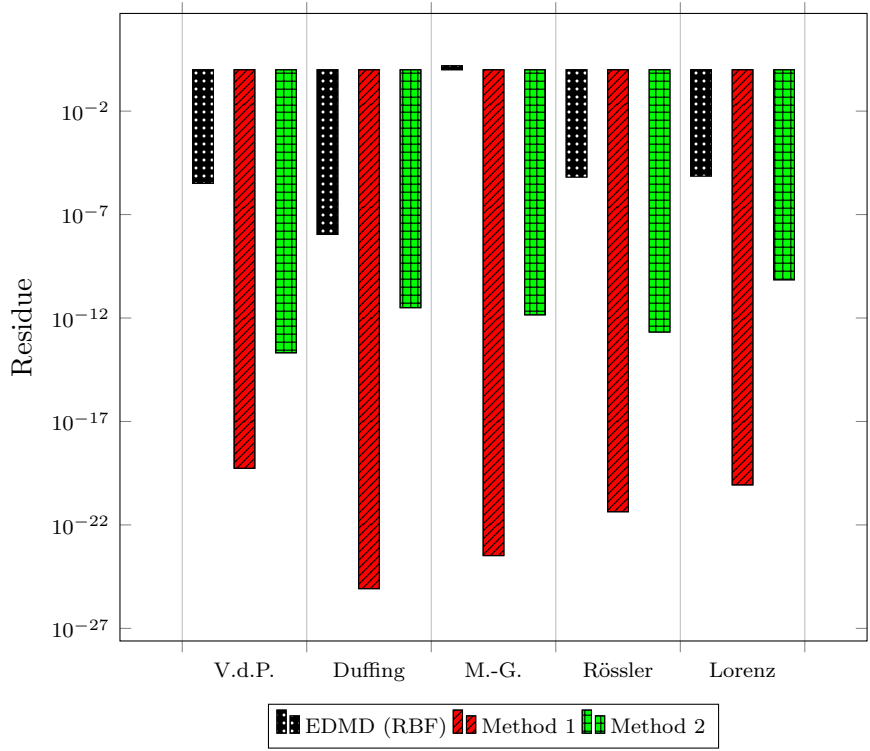

            The Koopman matrix computed with the basis functions generated by the reservoir is efficient to reconstruct the trajectories.
            As shown in FIG. \ref{fig:residues}, small residues \eqref{eq:residue} are obtained with all the systems introduced above.
            Moreover, in each case, the proposed methods outperform the EDMD method using as many radial basis functions as there are internal states in the reservoir.
            We observe that the first method yields better results than the second method.
            This can be explained by the fact the second method uses a smaller number $L<N$ of dictionary functions.
            We also note that the EDMD method is not able to reconstruct the trajectories of the Mackey-Glass system (residue larger than $1$).
            This suggests that the classical EDMD method with Gaussian radial basis functions cannot capture a time-delayed dynamics, in contrast to the reservoir whose internal states can be seen as functions depending on delayed input values (see Section \ref{subsec:meth1}).
            
\begin{figure}[H]
                \centering
                \begin{subfigure}{0.75\textwidth}
            		\centering
            		\includegraphics[width=0.9\columnwidth]{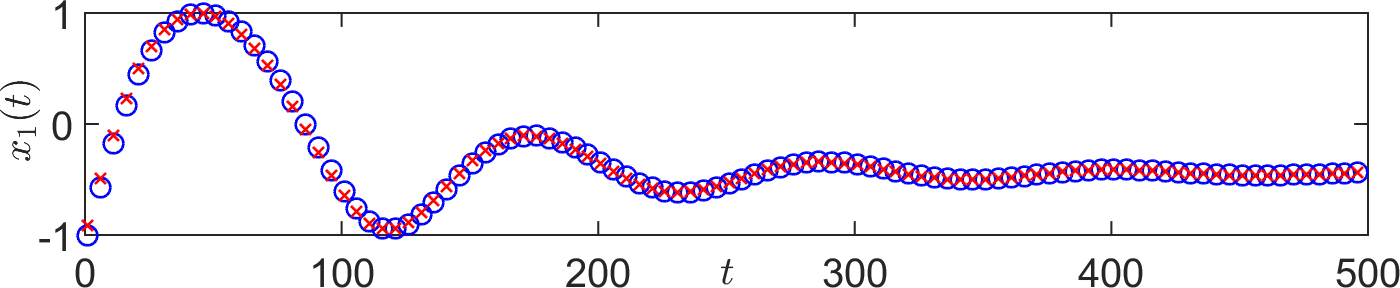}
            		\caption{}
            	\end{subfigure}
            	\begin{subfigure}{0.75\textwidth}
            		\centering
            		\includegraphics[width=0.9\columnwidth]{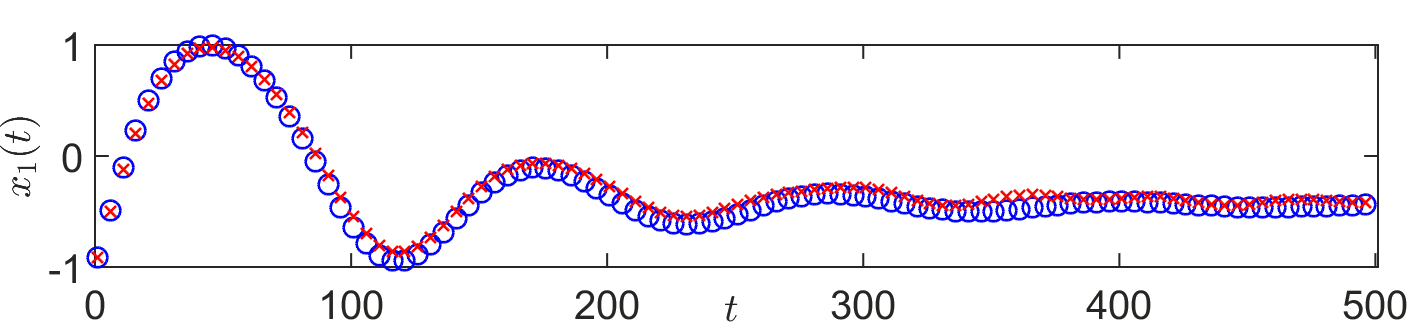}
            		\caption{}
            	\end{subfigure}	
            	\begin{subfigure}{0.75\textwidth}
            		\centering
            		\includegraphics[width=0.9\columnwidth]{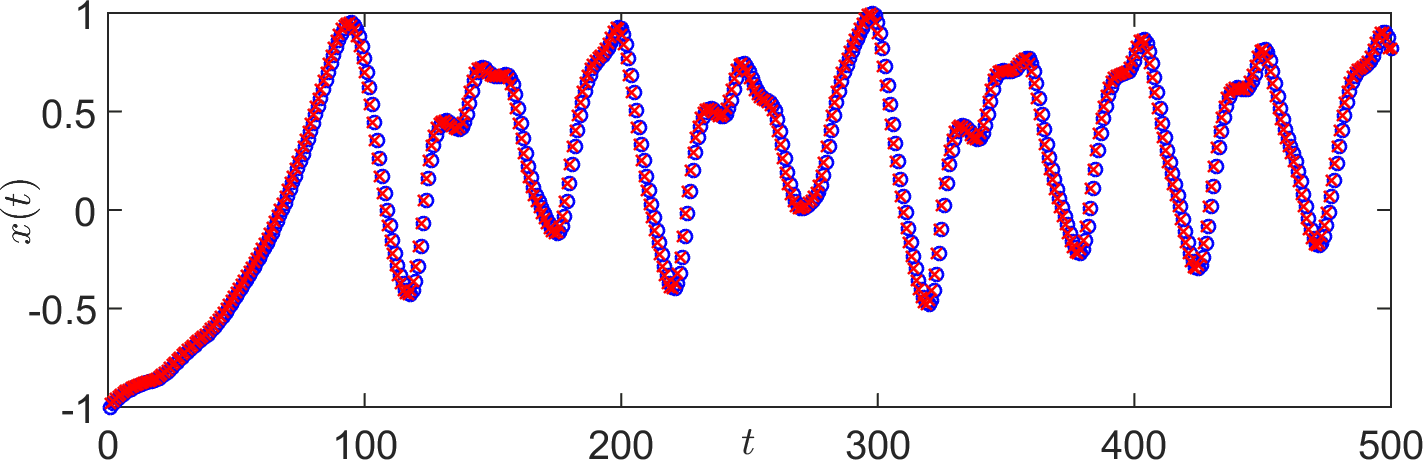}
            		\caption{}
            	\end{subfigure}	
            	\begin{subfigure}{0.75\textwidth}
            		\centering
            		\includegraphics[width=0.9\columnwidth]{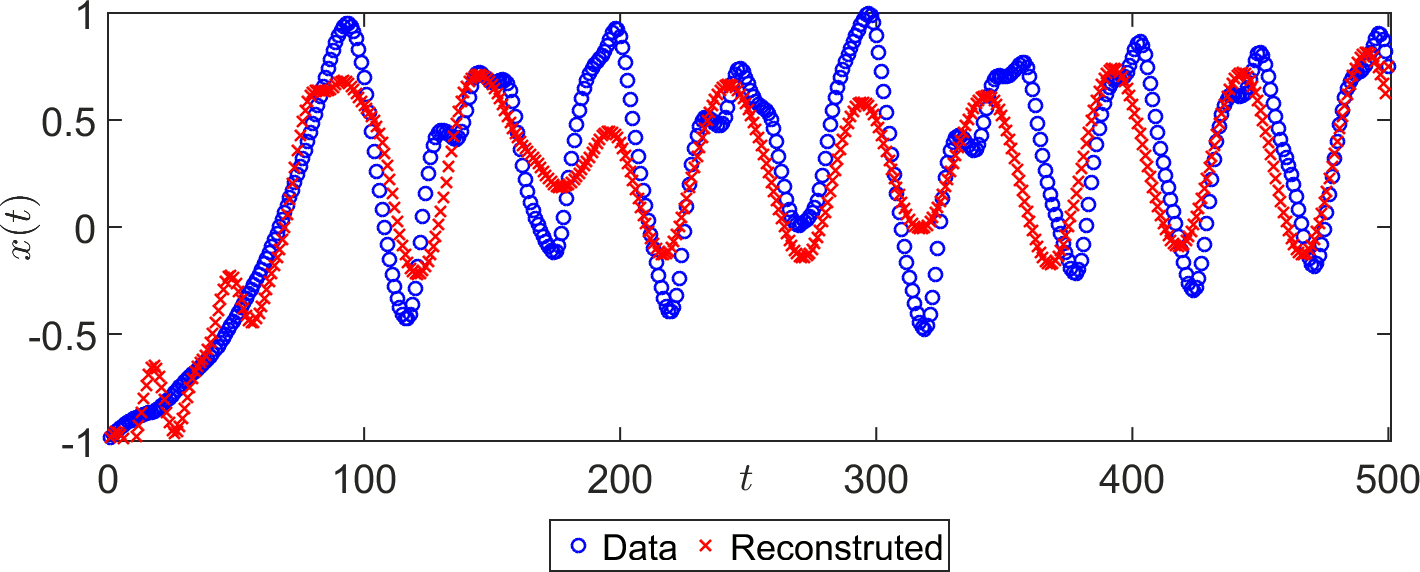}
            		\caption{}
            	\end{subfigure}	
            		\caption{\label{fig:recoBigFig}
            		The trajectories are reconstructed by iterating the Koopman matrix computed with Method 1 or 2 from the initial state. Blue circles and red crosses denote the data and the reconstructed trajectory, respectively (note that the all data points are not shown, so that the sampling period is smaller than it may seem on the figure).
            		The first two panels show the first component of the Duffing system reconstructed by Method 1 (panel (a)) and Method 2 (panel (b)).
            		The second component is not shown but is similar.
            		The last two panels show the trajectory of the Mackey-Glass system reconstructed by Method 1 (panel (c)) and Method 2 (panel (d)).
            		}
        	\end{figure}

            The proposed methods provide a Koopman matrix that can be iterated to reconstruct the trajectories from the initial state.
            This is illustrated by the $\bvec{nrmse}$ value shown in TABLE \ref{tab:recoNRMSEsummary}.
            We note that a larger error is obtained with the second method for the R\"{o}ssler system and the Lorenz system, in which case reconstructed trajectories diverge after some time.
            However, in all cases, the proposed methods outperform the EDMD method used with Gaussian radial basis functions.
            Finally, FIG. \ref{fig:recoBigFig} illustrates the reconstruction performances of the two methods for the Duffing system and the Mackey-Glass system.
            
            
            \begin{table}
                \centering
                \begin{tabular}{r|c|c|c|c|c}
                     \textbf{System}        & \textbf{V.d.P. }   & \textbf{Duffing}   & \textbf{M.-G. }    & \textbf{R\"{o}ssler}   & \textbf{Lorenz} \\
                     \hline
                     \textbf{EDMD (RBF)}    & $7.39 \times 10^{1}$      & $9.88 \times 10^{-7}$     & $3.43 \times 10^{0}$  & $9.32 \times 10^{-1}$ & $9.02 \times 10^{6}$ \\
                     \textbf{Method 1}      & $8.21 \times 10^{-11}$    & $5.32 \times 10^{-12}$    & $5.36 \times 10^{-8}$ & $2.80 \times 10^{-8}$ & $3.40 \times 10^{-9}$ \\
                     \textbf{Method 2}      & $2.37 \times 10^{-1}$     & $2.68 \times 10^{-1}$     & $1.84 \times 10^{-1}$ & $1.06 \times 10^{0}$  & $1.17 \times 10^{0}$ \\
                \end{tabular}
                \caption{\label{tab:recoNRMSEsummary}
                The mean value of the $\bvec{nrmse}$ vector components is shown for the different systems. The error vectors are computed according to \eqref{eq:error} for the first 100 reconstructed points. Our methods 1 and 2 yield better performance.}
            \end{table}
            
        \subsection{\label{subsec:prediciton}Prediction results}
        
            In this section, we briefly present prediction results for illustrative purposes. To do so, we consider the last data point of the training stage and iterate the Koopman matrix from this point.
            
\begin{figure}
            		\centering
            		\includegraphics[width=0.65\columnwidth]{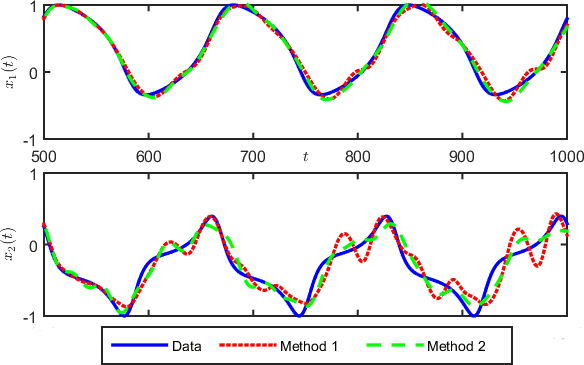}
            		\caption{\label{fig:predVdP} Both methods correctly predict the trajectory of the Van der Pol oscillator.}
            	\end{figure}

            \paragraph{Van der Pol system.}
            As a first toy example, we consider the Van der Pol oscillator and verify that both methods correctly predict the trajectory (FIG. \ref{fig:predVdP}). 
            For the second state variable, the prediction slowly diverges as the number of iterations of the Koopman matrix increases.
            Although the dynamics are very simple, it is noticeable that the classical EDMD method fails to provide good prediction results. In fact, for all study cases, the predicted trajectory either quickly diverges or converges to a constant.
            For this reason, we will not show the results obtained with the EDMD method.
         
\begin{figure}
            		\centering
            		\includegraphics[width=0.65\columnwidth]{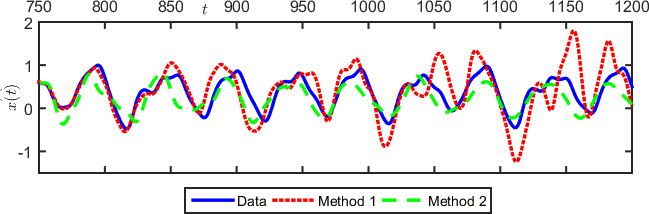}
            		\caption{\label{fig:predMG} The two methods are efficient to predict the trajectory of the Mackey-Glass dynamics over some time horizon.}
            	\end{figure}

\begin{figure}
            		\centering
            		\includegraphics[width=0.65\columnwidth]{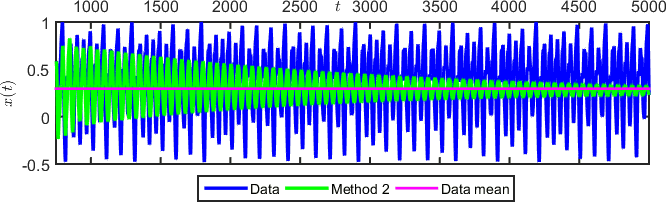}
            		\caption{\label{fig:pred2LongMG} The long-term prediction of the second method for the Mackey-Glass system converges to a constant value close to the mean value of the data.}
            	\end{figure}

            \paragraph{Mackey-Glass system.} FIG. \ref{fig:predMG} shows that the proposed methods are also efficient to predict the trajectory of the chaotic Mackey-Glass dynamics. 
            In the present setting, the first method provides more accurate results for short-term prediction, but the predicted trajectory eventually diverges. 
            In contrast the trajectory predicted with the second method converges to a constant which is close to the mean value of the data (see FIG. \ref{fig:pred2LongMG}).
            This is due to the fact that the Koopman matrix does not have unstable eigenvalues and has only the eigenvalue $\lambda = 1$ on the unit circle (see Section \ref{subsec:spectralProperties} below).
            It follows that the predicted trajectory converges to a steady state (associated with that eigenvalue $\lambda = 1$) which corresponds to the average value of the identity observable on the attractor, computed with respect to the stationary density.
            This value is also equal to the time average of the identity observable along a trajectory, which is close to the mean value of the data.
            
            	We finally recall that prediction is not the main goal of our methods.
            	Although the prediction results are decent, successive iterations of the Koopman matrix may lead to divergent prediction errors and could be avoided to improve the prediction results (see Section \ref{subsec:discussion}).

        \subsection{\label{subsec:spectralProperties}Spectral properties}
            
\begin{figure}
                \begin{subfigure}{0.4\textwidth}
            		\centering
            		\includegraphics[width=0.9\columnwidth]{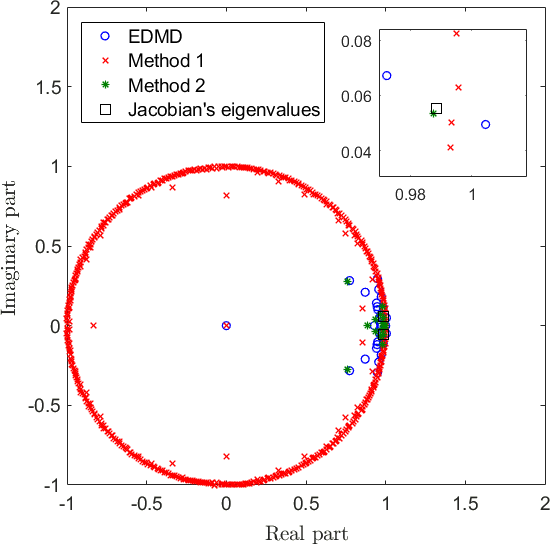}
            		\caption{\label{fig:spectrumDuffing}Eigenvalues of the Koopman matrix associated with the Duffing system.
            		}
            		\vspace{.75cm}
            	\end{subfigure}	
                \begin{subfigure}{0.4\textwidth}
            		\centering
            		\includegraphics[width=0.9\columnwidth]{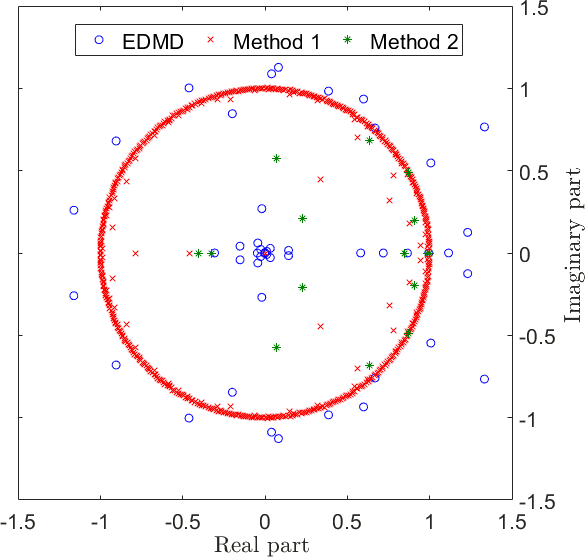}
            		\caption{\label{fig:spectrumRossler}
            		Zoom on the eigenvalues of the Koopman matrix associated with the R\"{o}ssler system.
            		}
            	\end{subfigure}
            	\caption{Computation of the spectrum of the Koopman operator for the Duffing system and the R\"{o}ssler system.
            }
        	\end{figure}

            In this Section, we compute an approximation of the spectrum of the Koopman operator for the Duffing system and the R\"{o}ssler system.
            The results are shown in FIG. \ref{fig:spectrumDuffing} and \ref{fig:spectrumRossler}, respectively.
            In both cases, the EDMD method generates many eigenvalues at the origin due to the rank-deficiency of the Koopman matrix.
            Our methods are characterized by less redundancy in the dictionary functions and, in particular, the second method provides a full-rank matrix.
            For the R\"{o}ssler system, the EDMD method also generates spurious eigenvalues outside the unit circle (eigenvalues with module approimately equal to $10$, not shown in the figure).
            This explains the fast divergence of the trajectories predicted with the EDMD method.
            In contrast, the first method recovers the whole unit circle, with a few additional eigenvalues inside the circle for both systems.
            The second method yields eigenvalues around $1$ and inside the unit circle for the Duffing system.
            It yields more scattered eigenvalues inside the unit circle for the R\"{o}ssler system.
            The inset in FIG. \ref{fig:spectrumDuffing} shows that the eigenvalues of the Jacobian matrix at the stable fixed points are  correctly recovered with the second method.
            %

        \subsection{\label{subsec:discussion}Discussion}
        
            \subsubsection{Comparison of the two methods}

            The first method provides the best reconstruction results since it exploits all the internal states of the reservoir rather than a few linear combinations of them.
            It should therefore be preferred if one aims at obtaining the most accurate matrix approximation of the Koopman operator.
            This method is also quite fast since there is no training process.
            
            However the Koopman matrix obtained with the first method may be very large, since the number of internal states of the reservoir computer is typically large. The second method is motivated by a tradeoff between the quality of the results and the size of the Koopman matrix at the cost of an additional computation time due to the dictionary training.
            It should be considered if one seeks for a low-dimensional approximation of the Koopman operator.
            In the context of prediction, it also seems that the second method provides slightly better results.
            While increasing the number of basis functions, the second method should produce results converging to those yielded by the first method.
            However the computation time also drastically increases so that the first method appears to be more relevant and efficient in this case. Similarly, if a large dataset is available, the second method might be computationally demanding because of huge matrices needed for the training.
                    
            \subsubsection{\label{subsubsec:strengthsWeaknesses}Strengths and weaknesses}
            
            A main advantage of the proposed methods is that they rely solely on linear techniques thanks to the reservoir computer framework, in contrast to other Koopman operator-based learning techniques.
            Moreover, both reconstruction and prediction results obtained with these methods are improved with respect to the results obtained with the classical EDMD method.
            The Koopman spectrum computed with these methods is consistent and motivates the use of a trained set of basis functions.
            Finally we note that both methods require very little data to provide accurate results.
            
            A main limitation of the second method is the computational cost of the dictionary training through the reservoir computer framework.
            We also note that both methods are not designed for prediction and cannot outperform state-of-the-art prediction methods.
            In particular the method fails to predict trajectories from initial conditions that are not related to the training set.
  
            \subsubsection{Improving the prediction methods}
            
            Our proposed methods mainly aim at computing the Koopman matrix with appropriate dictionary functions to provide the best global linear approximation of the dynamics.
            In the context of prediction, however, better results could be obtained with nonlinear approximations of the dynamics.
            We refer to other works proposing an efficient use of the reservoir computer for prediction, i.e. \cite{Pathak2018, gulina2018}.
                
            In the context of prediction with nonlinear approximations, we can also note that classical EDMD method could be used to compute a single iteration of the Koopman matrix, extract the updated value of the projection maps $\bvec{x}$, and use them to evaluate the iterated values of all dictionary functions. This would allow to project back the predicted trajectory onto the manifold containing the lifted states (see also a similar idea in the work \cite{bruder2019}).
            In fact, this amounts at computing the least squares projection of the system map $\bvec{F}$ in the span of the dictionary functions.
            Numerical simulations suggest that this method is very efficient in the context of prediction.
                
                
	\section{\label{sec:Conclusions}Conclusion}
	
    	We have proposed two novel methods for computing a finite-dimensional approximation of the Koopman operator.
    	These methods combine classical EDMD with the use of a reservoir computer.
    	In the first method, the dictionary functions are chosen to be the internal states of the reservoir.
    	In the second method, the reservoir computer is trained and the dictionary functions are optimized linear combinations of internal states.
    	A key advantage of these two methods is that they rely on linear optimization techniques.
    	The accuracy of the Koopman matrix approximation is assessed in the context of reconstruction, prediction, and computation of the Koopman operator spectrum.
    	The results are encouraging and pave the way to the use of the reservoir computer in the Koopman operator framework.
    	
        Several research perspectives can be proposed.
        First, the method could be improved to achieve better predictive performances, although this is not our main goal in this paper (note also that other Koopman operator-based methods exist for this purpose, see e.g. \cite{giannakis2018,Korda2018,Khodkar2019}).
        To do so, one could adapt the training of the reservoir according to this specific prediction objective, promote the computation of stable Koopman matrices, and use proper projections between the iterations of the Koopman matrix (see e.g. \cite{bruder2019}).
        Our second method could also be complemented with convergence results for the alternating optimization scheme.
        Finally, the proposed methods could be used on real datasets, in the context of spectral analysis, network identification, time-series classification, event detection, and predictive control.
    
	\section{\label{sec:Acknowledgments}Acknowledgments}
	    This work is supported by the Namur Institute For Complex Systems (naXys) at University of Namur.
	    M. Gulina is grateful to Adrien Fiorucci and Thomas Evrard for fruitful discussions and comments on the manuscript.
	   
	\section*{Data availability}
		The data that support the findings of this study are almost all available within the article.
		Any data that are not available can be found from the corresponding author upon reasonable request.

    \bibliographystyle{siam}
    \bibliography{chaos2020}

\begin{thebibliography}{10}

\bibitem{gulina2018}
{\sc P.~Antonik, M.~Gulina, J.~Pauwels, and S.~Massar}, {\em Using a reservoir
  computer to learn chaotic attractors, with applications to chaos
  synchronization and cryptography}, Phys. Rev. E, 98 (2018), p.~012215.

\bibitem{ArabiMezic2017}
{\sc H.~Arbabi and I.~Mezić}, {\em {Ergodic Theory, Dynamic Mode
  Decomposition, and Computation of Spectral Properties of the Koopman
  Operator}}, SIAM Journal on Applied Dynamical Systems, 16 (2017),
  p.~2096–2126.

\bibitem{bollt2020}
{\sc E.~Bollt}, {\em {On Explaining the surprising success of reservoir
  computing forecaster of chaos? The universal machine learning dynamical
  system with contrasts to VAR and DMD}}, arXiv preprint arXiv:2008.06530,
  (2020).

\bibitem{bruder2019}
{\sc D.~Bruder, B.~Gillespie, C.~D. Remy, and R.~Vasudevan}, {\em Modeling and
  control of soft robots using the koopman operator and model predictive
  control}, arXiv preprint arXiv:1902.02827,  (2019).

\bibitem{dogra2020optimizing}
{\sc A.~S. Dogra and W.~T. Redman}, {\em {Optimizing neural networks via
  Koopman operator theory}}, arXiv preprint arXiv:2006.02361,  (2020).

\bibitem{Drmac2017}
{\sc Z.~Drmač, I.~Mezic, and R.~Mohr}, {\em Data driven modal decompositions:
  Analysis and enhancements}, SIAM Journal on Scientific Computing, 40 (2017).

\bibitem{giannakis2018}
{\sc D.~Giannakis, S.~Das, and J.~Slawinska}, {\em {Reproducing kernel Hilbert
  space compactification of unitary evolution groups}}, arXiv preprint
  arXiv:1808.01515,  (2018).

\bibitem{H_Tu_2014}
{\sc J.~H.~Tu, C.~W.~Rowley, D.~M.~Luchtenburg, S.~L.~Brunton, and
  J.~Nathan~Kutz}, {\em On dynamic mode decomposition: Theory and
  applications}, Journal of Computational Dynamics, 1 (2014), p.~391–421.

\bibitem{Huang2018}
{\sc B.~{Huang}, X.~{Ma}, and U.~{Vaidya}}, {\em Feedback stabilization using
  koopman operator}, in 2018 IEEE Conference on Decision and Control (CDC),
  2018, pp.~6434--6439.

\bibitem{jaeger2001echo}
{\sc H.~Jaeger}, {\em The “echo state” approach to analysing and training
  recurrent neural networks}, GMD-Report 148, German National Research
  Institute for Computer Science,  (2001).

\bibitem{jaeger2004harnessing}
{\sc H.~Jaeger and H.~Haas}, {\em Harnessing nonlinearity: Predicting chaotic
  systems and saving energy in wireless communication}, Science, 304 (2004),
  pp.~78--80.

\bibitem{kaiser2017}
{\sc E.~Kaiser, J.~N. Kutz, and S.~L. Brunton}, {\em {Data-driven discovery of
  Koopman eigenfunctions for control}}, arXiv preprint arXiv:1707.01146,
  (2017).

\bibitem{Khodkar2019}
{\sc M.~A. Khodkar, P.~Hassanzadeh, and A.~Antoulas}, {\em {A Koopman-based
  framework for forecasting the spatiotemporal evolution of chaotic dynamics
  with nonlinearities modeled as exogenous forcings}}, arXiv preprint
  arXiv:1909.00076,  (2019).

\bibitem{Koopman1931}
{\sc B.~O. Koopman}, {\em {Hamiltonian Systems and Transformation in Hilbert
  Space}}, Proceedings of the National Academy of Sciences of the United States
  of America, 17 (1931), pp.~315--318.
\newblock 16577368[pmid].

\bibitem{Korda2018}
{\sc M.~Korda and I.~Mezić}, {\em {Linear predictors for nonlinear dynamical
  systems: Koopman operator meets model predictive control}}, Automatica, 93
  (2018), pp.~149 -- 160.

\bibitem{Kevrekidis2017}
{\sc Q.~Li, F.~Dietrich, E.~M. Bollt, and I.~G. Kevrekidis}, {\em {Extended
  dynamic mode decomposition with dictionary learning: A data-driven adaptive
  spectral decomposition of the Koopman operator}}, Chaos: An Interdisciplinary
  Journal of Nonlinear Science, 27 (2017), p.~103111.

\bibitem{lukovsevivcius2012reservoir}
{\sc M.~Luko{\v{s}}evi{\v{c}}ius, H.~Jaeger, and B.~Schrauwen}, {\em Reservoir
  computing trends}, KI - K{\"u}nstliche Intelligenz, 26 (2012), pp.~365--371.

\bibitem{lukovsevivcius2009survey}
{\sc M.~Luko\v{s}evi\v{c}ius and H.~Jaeger}, {\em Survey: Reservoir computing
  approaches to recurrent neural network training}, Comput. Sci. Rev., 3
  (2009), p.~127–149.

\bibitem{Lusch2018}
{\sc B.~Lusch, J.~N. Kutz, and S.~L. Brunton}, {\em Deep learning for universal
  linear embeddings of nonlinear dynamics}, Nature Communications, 9 (2018),
  p.~4950.

\bibitem{MackeyGlassSystem}
{\sc M.~Mackey and L.~Glass}, {\em Oscillation and chaos in physiological
  control systems}, Science (New York, N.Y.), 197 (1977), p.~287—289.

\bibitem{Manojlovic2020}
{\sc I.~Manojlović, M.~Fonoberova, R.~Mohr, A.~Andrejčuk, Z.~Drmač,
  Y.~Kevrekidis, and I.~Mezić}, {\em {Applications of Koopman Mode Analysis to
  Neural Networks}}, arXiv preprint arXiv:2006.11765,  (2020).

\bibitem{Mauroy2017}
{\sc A.~Mauroy and J.~Hendrickx}, {\em Spectral identification of networks
  using sparse measurements}, SIAM Journal on Applied Dynamical Systems, 16
  (2017), pp.~479--513.
\newblock Article.

\bibitem{Mauroy2016GlobalStability}
{\sc A.~{Mauroy} and I.~{Mezić}}, {\em {Global stability analysis Using the
  eigenfunctions of the Koopman operator}}, IEEE Transactions on Automatic
  Control, 61 (2016), pp.~3356--3369.

\bibitem{Mauroy2013}
{\sc A.~Mauroy, I.~Mezi{\'{c}}, and J.~Moehlis}, {\em {Isostables, isochrons,
  and Koopman spectrum for the action--angle representation of stable fixed
  point dynamics}}, Physica D: Nonlinear Phenomena, 261 (2013), pp.~19--30.

\bibitem{mauroy2020koopman}
{\sc A.~Mauroy, y.~Susuki, and I.~Mezi{\'{c}}}, {\em The Koopman Operator in
  Systems and Control}, Springer, 2020.

\bibitem{Mezic2005}
{\sc I.~Mezi{\'{c}}}, {\em Spectral properties of dynamical systems, model
  reduction and decompositions}, Nonlinear Dynamics, 41 (2005), pp.~309--325.

\bibitem{Pan2020}
{\sc S.~Pan and K.~Duraisamy}, {\em Physics-informed probabilistic learning of
  linear embeddings of nonlinear dynamics with guaranteed stability}, SIAM
  Journal on Applied Dynamical Systems, 19 (2020), p.~480–509.

\bibitem{Pathak2018}
{\sc J.~Pathak, A.~Wikner, R.~Fussell, S.~Chandra, B.~R. Hunt, M.~Girvan, and
  E.~Ott}, {\em Hybrid forecasting of chaotic processes: Using machine learning
  in conjunction with a knowledge-based model}, Chaos: An Interdisciplinary
  Journal of Nonlinear Science, 28 (2018), p.~041101.

\bibitem{Proctor2015}
{\sc J.~L. Proctor and P.~A. Eckhoff}, {\em Discovering dynamic patterns from
  infectious disease data using dynamic mode decomposition}, International
  health, 7 (2015), pp.~139--145.
\newblock 25733564[pmid].

\bibitem{Rowley2009}
{\sc C.~Rowley, I.~Mezic, S.~Bagheri, P.~Schlatter, and D.~Henningson}, {\em
  Spectral analysis of nonlinear flows}, Journal of Fluid Mechanics, 641
  (2009), pp.~115 -- 127.

\bibitem{Schmid2010}
{\sc P.~J. Schmid}, {\em {Dynamic mode decomposition of numerical and
  experimental data}}, Journal of Fluid Mechanics, 656 (2010), pp.~5--28.

\bibitem{Skibinsky-Gitlin2018}
{\sc E.~Skibinsky-Gitlin, M.~Alomar, C.~Frasser, V.~Canals, E.~Isern, M.~Roca,
  and J.~L. Rossello}, {\em Cyclic Reservoir Computing with FPGA Devices for
  Efficient Channel Equalization}, 01 2018, pp.~226--234.

\bibitem{Surana2018}
{\sc A.~Surana}, {\em {Koopman operator framework for time series modeling and
  analysis}}, Journal of Nonlinear Science,  (2018).

\bibitem{Susuki2014}
{\sc Y.~{Susuki} and I.~{Mezić}}, {\em {Nonlinear Koopman modes and power
  system stability assessment without models}}, IEEE Transactions on Power
  Systems, 29 (2014), pp.~899--907.

\bibitem{SusukiMezic2015}
\leavevmode\vrule height 2pt depth -1.6pt width 23pt, {\em {A prony
  approximation of Koopman Mode Decomposition}}, in 2015 54th IEEE Conference
  on Decision and Control (CDC), 2015, pp.~7022--7027.

\bibitem{Susuki2011}
{\sc Y.~{Susuki} and I.~{Mezi{\'{c}}}}, {\em {Nonlinear Koopman modes and
  coherency identification of coupled swing dynamics}}, IEEE Transactions on
  Power Systems, 26 (2011), pp.~1894--1904.

\bibitem{Takeishi2017}
{\sc N.~Takeishi, Y.~Kawahara, and T.~Yairi}, {\em {Learning Koopman invariant
  subspaces for Dynamic Mode Decomposition}}, Advances in Neural Information
  Processing Systems (Proc. of NIPS), 30 (2017).

\bibitem{Triefenbach2010}
{\sc F.~Triefenbach, A.~Jalalvand, B.~Schrauwen, and J.-P. Martens}, {\em
  Phoneme recognition with large hierarchical reservoirs}, in Advances in
  Neural Information Processing Systems, J.~Lafferty, C.~Williams,
  J.~Shawe-Taylor, R.~Zemel, and A.~Culotta, eds., vol.~23, Neural Information
  Processing System Foundation, 2010, p.~9.

\bibitem{Williams2015}
{\sc M.~O. Williams, I.~G. Kevrekidis, and C.~W. Rowley}, {\em {A data--driven
  approximation of the Koopman operator: Extending Dynamic Mode
  Decomposition}}, Journal of Nonlinear Science, 25 (2015), pp.~1307--1346.

\bibitem{Yeung2019}
{\sc E.~{Yeung}, S.~{Kundu}, and N.~{Hodas}}, {\em {Learning deep neural
  network representations for Koopman operators of nonlinear dynamical
  systems}}, in 2019 American Control Conference (ACC), 2019, pp.~4832--4839.

\end{thebibliography}
	
\end{document}